\documentclass[final,leqno]{siamltex}

\usepackage{amsfonts,amsmath,amssymb}
\usepackage{caption}
\usepackage{color}
\usepackage{subcaption,tikz}
\usetikzlibrary{matrix}
\def\card{{card}}
\newtheorem{remark}[theorem]{Remark}

\title{Ordered Field Property for Zero-Sum Stochastic Games\thanks{This research was supported by The Australian Research Council Discovery Grant DP150100618. The third author wishes to thank T.E.S. Raghavan for introducing him to this problem in the late $1970's$.}}

\author{K. Avrachenkov\footnotemark[2]\footnotemark[2]
\and V. Ejov\footnotemark[3]
				\and J. A. FILAR\footnotemark[5],\footnotemark[6]
				\and A. Moghaddam\footnotemark[5],\footnotemark[7]}

\begin{document}
\maketitle

\renewcommand{\thefootnote}{\fnsymbol{footnote}}

\footnotetext[2]{Inria Sophia Antipolis,\ {\tt k.avrachenkov@inria.fr}}
\footnotetext[3]{Collage of Science and Engineering Flinders University of South Australia
				 Bedford Park SA5042\newline
MSU, Faculty of Mechanics and Mathematics, Russia,119991,
				 Moscow, GSP-1, 1 Leninskiye Gory}
\footnotetext[5]{Centre for Applications in Natural Resource Mathematics,
								School of Mathematics and Physics,
								The University of Queensland,
								St Lucia, Queensland, 4072}
\footnotetext[6]{{\tt j.filar@uq.edu.au}}
\footnotetext[7]{{\tt a.moghaddam1@uq.edu.au}}
%\footnotetext[5]{Address of A.~U. Thorthree}
%\footnotetext[6]{Support in common for the first and second
%authors.}

\renewcommand{\thefootnote}{\arabic{footnote}}

\begin{abstract}
We consider a finite state, finite action, zero-sum stochastic games with data defining the game lying in the ordered field of algebraic numbers. In both the discounted and the limiting average versions of these games we prove that the value vector also lies in the same field of algebraic numbers. In a prescribed sense, our results settle a problem that has remained open since, at least, $1991$.
\end{abstract}

\begin{keywords}Stochastic games, ordered field property, algebraic numbers, algebraic variety, Gr{\" o}bner basis, polynomial equations. \end{keywords}

\begin{AMS} 90D15 \end{AMS}

\pagestyle{myheadings}
\thispagestyle{plain}
\markboth{Ordered field property and stochastic games}{Ordered field property and stochastic games}

 %*****************************************************************************************************************************%
 %**************************************************--- Introduction ---*******************************************************%
 %*****************************************************************************************************************************%

\section{Introduction}

Arguably, modern game theory was launched in $1928$ by von Neumann (see \cite{Neun28}). His seminal paper proved that a finite, zero-sum, two person (matrix) game possesses a game-theoretic value and a pair of optimal strategies. Subsequently, Weyl (see \cite{Weyl}) supplied a simpler proof and, in addition, showed that if the entries of such a matrix game belong to an ordered field, then the value belongs to the same ordered field.

Stochastic games were introduced in 1953 by Shapley (see \cite{Shapley}). Shapley's formulation was analogous to what are now called discounted stochastic games. In these games the rewards at future stages are discounted by a factor $\beta \in [0, 1)$. However, in a remark, Shapley also observed that Weyl's ordered field property does not hold for stochastic games. In particular, when the data of these games lie in the field of rational numbers, the value vector need not be in the same field.

The latter remark, stimulated a whole line of research attempting to characterize special classes of stochastic games that possess the ordered field property.  Indeed, at least five, structured, classes of stochastic games have been identified and shown to possess the ordered field property, over the field $\mathbb{Q}$ of rational numbers.  These classes include stochastic games of perfect information originally introduced by Gillette in 1957 \cite{Gillette1957}, separable reward and state independent transition (SER-SIT) games \cite{Sobel1981myopic}, \cite{Parthasarathy1984}, single-controller games \cite{Parthasarathy1981}, switching-controller games \cite{Filar1981} and additive reward and additive transition (ARAT) games \cite{Raghavan85}.  However, a more comprehensive characterization of the ordered field property proved to be challenging and was named as one of the open problems in the topic in a 1991 survey paper \cite{Raghavan1991}.

In this paper we provide a rather complete answer to the above mentioned open problem, in the sense that we, constructively, establish the ordered field property over the ordered field $\mathbb{F}$ of real algebraic numbers. This is shown without any structural assumptions on the game and applies to both discounted and limiting average stochastic games that will be introduced later on.  Arguably, $\mathbb{F}$ is the most natural ordered field to consider in view of the fact that the field of rational numbers $\mathbb{Q}$ was too small.  This is because $\mathbb{F}$ is a finite degree field extension of $\mathbb{Q}$ is also countable and, in the sense of Rotman (\cite{rotman1998}, page 50), it is smallest field containing rationals. Our analysis exploits Gr{\" o}bner basis methods as well as results reported in Szczechla et al. \cite{Filar97} and we also use similar notation to the latter paper.

We note that the asymptotic behavior of the value vector as $\beta\to 1$ has been studied by a number of authors. In an important, and related,  contribution Bewley and Kohlberg \cite{BeKo} viewed Shapley's ``optimality condition'' as an elementary sentence in formal logic over the closed ordered field of real Puiseux series. These authors invoke a powerful theorem from mathematical logic, known as Tarski's principle, to conclude that in some neighborhood of $\beta = 1$, the value vector belongs to the field of real Puiseux series. We also note that our approach has similarities to the semialgebraic theory of stochastic games proposed in Milman \cite{Mil2002} and continued by Neyman \cite{Ney2003}. The latter permits analysis of sets defined by both polynomial equations and inequalities but it is not immediately clear how to apply Gr{\"o}bner basis approach in that context.

%The essentially algebraic nature of Bewley and Kohlberg’s approach and the use of Tarski’s principle, while ingenious, do not give insight into the manner in which fractional power series solutions arise naturally in this problem. Arguably, this has proved to be a difficulty for researchers in stochastic games, because Bewley and Kohlberg’s result has become an important building block in subsequent developments. For instance, Mertens and Neyman \cite{Mertens} used it to prove the existence of the value vector in the limiting average stochastic games.

\section{Definitions and preliminaries of matrix games}

Any $m\times n$ real matrix $A = (a_{ij})^{m,n}_{i,j=1}$ can be regarded as a two-person, zero-sum matrix game with $a_{ij}$ denoting the amount player $II$ will pay player $I$ if $II$ chooses an action $j\in{1, 2,\cdots, n}$ and $I$ chooses an action $i\in{1, 2,\cdots, m}$. A mixed (or randomized) strategy for player $I (II)$ in such a game is an m(n)-component probability vector x(y) whose $i^{th} (j^{th})$ entry $x_i(y_j)$ denotes the probability that player $I (II)$ will choose an action $i(j)$. A consequence of the celebrated ``minimax theorem'' for matrix games \cite{Neun28} is  that there always exists a strategy pair $(x^{0}, y^{0})$ satisfying:
\begin{equation}
x^T A y^{0}\leq (x^{0})^T A y^{0}\leq (x^{0})^T A y
\end{equation}
for all mixed strategies $x(y)$ of player $I (II)$. The strategies $(x^{0}, y^{0})$ are then called \textit{optimal strategies}, and the real number $val(A) := (x^{0})^T A y^{0}$ is called the \textit{value} of the matrix game A. It is well known that if $b_{ij} = k a_{ij} + c$, for all $i, j, k > 0$, and $B = (b_{ij})^{m,n}_{i,j=1}$, then
$valB = k valA+c$. Hence, there is no loss of generality in assuming that the value of a matrix game is either positive, or simply, not equal to zero.

A matrix game A is called \textit{completely mixed} if all of its optimal strategies are strictly positive in every component. Extending earlier results of Kaplansky \cite{Kap}, Shapley and Snow \cite{ShapSnow} demonstrated the following result.

\begin{proposition}
If A is a matrix game and $val A \neq 0$, then $A$ has a square invertible submatrix $\bar{A}$, called a {\em Shapley--Snow kernel}, such that
\begin{enumerate}
	\item[(K1)] $val\bar{A}=val A=\frac{\det(\bar{A})}{\sum_{ij}\bar{A}^{ij}}=\left(\mathbf{1}^T[\bar{A}^{-1}]\mathbf{1}\right)^{-1}$ where $\bar{A}^{ij}$ denotes the $\left(i,j\right)^{th}$ cofactor of $\bar{A}$.
	\item[(K2)] There is a pair $\left(x^0, y^0\right)$ of strategies for $\bar{A}$ which is optimal for $\bar{A}$ (after inserting zeroes) and satisfies $\left(x^0\right)^T=\left(val A\right) \mathbf{1}^T[\bar{A}^{-1}]$ and $y^0=\left(val A\right)[\bar{A}^{-1}]\mathbf{1}$.
\end{enumerate}
\end{proposition}

\begin{remark}
If $val A \neq 0,$  Lemma $1.2$ of \cite{Filar97} shows that it is always possible to find a Shapley-Snow kernel $\bar{A}$ that is a completely mixed matrix game and, of course, still satisfies $(K1)-(K2),$ above. Following \cite{Filar97} we shall refer to such a kernel as {\em cmv-Shapley-Snow kernel,} or simply {\em cmv-kernel}. Their usefulness stems from the fact that the algebraic formula for their value $val\bar{A}=\left(\mathbf{1}^T[\bar{A}^{-1}]\mathbf{1}\right)^{-1}$ is invariant under small perturbations.
\end{remark}
\newpage
 %**************************************************************************************************************************%
 %**************************************************--- New Section ---*****************************************************%
 %**************************************************************************************************************************%

\section{Discounted stochastic games}
\subsection{Definitions and preliminaries of stochastic games}
A \textit{stochastic game} as formulated by Shapley \cite{Shapley} is played in stages. At each stage, the game is in one of finitely many states, $s=1,2,\ldots,N$, in which players $I$ and $II$ are obliged to play a matrix game $R(s)=(r(s,a,b))^{m_s,n_s}_{a,b=1}$, once. The ``law of motion" is defined by $p(s'|s,a,b)$, where the event $\left\{s'|s,a,b\right\}$ is the event that the game will enter state $s'$ at the next stage given that at the current stage the state of the game is $s$, and that players $I$ and $II$ choose the $a^{th}$ row and the $b^{th}$ column of $R(s)$, respectively.

In general, players' strategies will depend on complete past histories. In this paper, however, we shall only be concerned with \textit{stationary strategies}. We may represent a typical stationary strategy $\mu$ for player $I$ by a ``composite" vector, $\mu=(\mu(1),\mu(2),\ldots,\mu(N))$, where each $\mu(s)$ is a probability distribution on $\left\{1,2, \ldots, m_s\right\}$. Player $II$'s stationary strategies $\nu$ are similarly defined.

It should be clear that once we specify the initial state s and a strategy pair $(\mu,\nu)$ for players $I$ and $II$, we implicitly define a probability distribution over all sequences of states and actions which can occur during the game and consequently over all sequences of payoffs to player $I$. In particular, if the random variable $\mathcal{R}_t$ denotes the payoff to player $I$ at stage t, then the expected value of $\mathcal{R}_t$ given $s$ and $(\mu,\nu)$
\begin{equation}
E_{\mu\nu s}(\mathcal{R}_t):=E\lbrace\mathcal{R}_t|\mu,\nu,s\rbrace
\end{equation}
is well defined. The $\beta$\textit{-discounted stochastic game} $\Gamma_{\beta}$ is then the game in which the overall payoff, normalized by a factor of $1-\beta$, resulting from the strategy pair $(\mu,\nu)$ and a starting state s is evaluated according to
\begin{eqnarray}
v_{\beta}(\mu,\nu,s):= \sum_{t=1}^{\infty} \beta^{t-1}(1-\beta)E_{\mu\nu s}(\mathcal{R}_t), \nonumber
\end{eqnarray}
where $\beta \in (0,1)$ is called the discount factor. The number $v_s(\beta)$ is called the \textit{value} of the game $\Gamma_{\beta}$ starting in state $s$ if $v_s(\beta) = \sup_{\mu} \inf_{\nu}v_{\beta}(\mu,\nu,s) = \inf_{\nu} \sup_{\mu} v_{\beta}(\mu,\nu,s)$. The vector $\mathbf{v}(\beta) = (v_1(\beta), v_2(\beta),\ldots, v_N(\beta))$ is called the \textit{value vector}. Furthermore, the pair $(\mu^0, \nu^0)$ is called an optimal strategy pair for players $I$ and $II$ if, for each starting state $s$,
\begin{eqnarray}
v_s(\beta)=v_{\beta}(\mu^0, \nu^0,s).\nonumber
\end{eqnarray}
The existence of the value vector and of a pair of \textit{optimal stationary strategies} was proved in 1953 in Shapley's seminal paper on the subject \cite{Shapley}. A key element in Shapley's proof was the construction of $N$ auxiliary matrix games $R_{\beta}(s,\mathbf{u})$ that depend on an arbitrary vector
$\mathbf{u}=(u_1,u_2,\ldots,u_N)\in\mathbb{R}^\mathbf{N}$ according to
\begin{equation}
R_{\beta}(s,\mathbf{u})=\left[(1-\beta)r(s,a,b)+\beta\sum_{s'=1}^N p(s'|s,a,b)u_{s'}\right]_{a,b=1}^{m_s,n_s}.
\end{equation}
In view of the fact that the value of a matrix game always exists, it is possible to define, for each $\beta \in (0,1)$, an operator $T_{\beta}:\mathbb{R}^{\mathbf{N}}\longrightarrow \mathbb{R}^{\mathbf{N}}$, the $s^{th}$ component of which is given by
\begin{equation}
[T_{\beta}(\mathbf{u})]_s:=val[R_{\beta}(s,\mathbf{u})].
\end{equation}
This operator is a contraction operator in the sup-norm with contraction constant $\leq\beta$; see \cite{Shapley}. It therefore follows from Banach's fixed-point theorem that there exists a unique fixed point $\mathbf{v}(\beta)$ of $T_{\beta}$; that is,
\begin{equation}\label{FP}
\mathbf{v}(\beta)=T_{\beta}(\mathbf{v}(\beta)).
\end{equation}
%This can be shown to be an equivalent definition of the value vector $\mathbf{v}(\beta)$ introduced
%above.
Also, any set of optimal strategy pairs for $R_{\beta}(s,\mathbf{v}(\beta))$ ($s=1,2,\ldots,N$) can be shown to form an
optimal strategy pair for $\Gamma_{\beta}$.

%Since ${T_{\beta}}$ is a continuous family of contractions we have the following result.
%
%\begin{lemma}
%The function $\beta \mapsto\mathbf{v}(\beta)$ is bounded and continuous on $(0,1)$.
%\end{lemma}

Assume $v_s(\beta)\neq 0$ for all $\beta$ and $s$.  In view of Proposition 2.1 we know that for each fixed $\beta\in (0,1)$ and each $\mathbf{u}$ close enough to $\mathbf{v}(\beta)$ there exist cmv-kernels $\bar{R}_{\beta}(s,\mathbf{u})$ such that the fixed-point equation above reduces to
\begin{equation}\label{Redusedfixedpointeq}
v_s(\beta)=\frac{\det(\bar{R}_{\beta}(s,\mathbf{v}(\beta)))}{\sum_i\sum_j[\bar{R}_{\beta}(s,\mathbf{v}(\beta))]^{ij}}\ \ \textrm{for each } s=1,2,\ldots, N,
\end{equation}
where $[V]^{ij}$ denotes the $(i,j)^{th}$ cofactor of a matrix $V$.

If we now transform the above equations to
\begin{equation}\label{Redusedfixedpointeq2ndform}
v_s(\beta)\left\{\sum_i\sum_j[\bar{R}_{\beta}(s,\mathbf{v}(\beta))]_{ij}\right\}
-\det(\bar{R}_{\beta}(s,\mathbf{v}(\beta)))=0,\ s=1,2,\ldots, N,
\end{equation}
then for each fixed combination of the locations of these kernels we can regard this system as being a system of polynomials in the variables $z_0:=\beta, z_1:=v_1(\beta), z_2:=v_2(\beta), \ldots, z_N:=v_N(\beta)$ of the form

\begin{equation}\label{syseofcoupledpolyeeqes}
 		\left\{
 		\begin{array}{c}
 		f_1(z_0,\ldots,z_N)=0\\
		\vdots\ \\
		f_N(z_0,\ldots,z_N)=0.
 		\end{array}
 		\right.
\end{equation}

Consider any fixed collection of non-empty subsets of actions $K_s \subseteq \lbrace 1, \ldots, m_s\rbrace$ and $L_s \subseteq \lbrace 1, \ldots, n_s\rbrace$ such that $\card\ K_s = \card\ L_s$ for all $s = 1, \ldots, N$. Thus we have a finite sequence of sets $\kappa = (K_1, L_1,\ldots, K_N, L_N)$. By restricting the players' actions in each state $s$ to $K_s$ and $L_s$, we obtain a discounted stochastic game, $\Gamma_{\beta,\kappa}$ with auxiliary Shapley games
\begin{equation}
R_{\beta,\kappa}(s,\mathbf{u})=[R_{\beta}(s,\mathbf{u})_{ij}]_{i\in K_s, j\in L_s}.
\end{equation}
For fixed $\kappa$, let
\begin{equation}
\bar{R}_{\beta}(s,\mathbf{u})=R_{\beta,\kappa}(s,\mathbf{u})=[R_{\beta}(s,\mathbf{u})_{ij}]_{i\in K_s, j\in L_s}.
\end{equation}
be the $s^{th}$ cmv-kernel associated with $\kappa$. Clearly the system of polynomials (\ref{syseofcoupledpolyeeqes}) depends on $\kappa$ but this dependence is suppressed in order to simplify already complicated notation. However, it is important to note that $K=\lbrace \kappa | \text{defining cmv-kernels}\rbrace$ is finite. Hence, there are at most $|K|$ systems of polynomials of the form (\ref{syseofcoupledpolyeeqes}) to consider, where $|K|$ denotes the cardinality of $K$.

 %**************************************************************************************************************************%
 %**************************************************--- New Section ---*****************************************************%
 %**************************************************************************************************************************%

\subsection{Discounted stochastic games over the ordered field of algebraic numbers}

Suppose now that all the data of the discounted stochastic game $\Gamma_\beta$, namely, $r(s,a,b)$, $p(s'|s,a,b)$ (for all, $s, s', a, b$) and $\beta$ lie in the field $\mathbb{F}$ of real algebraic numbers.
From the outset, we wish to outline a very elegant - but non-constructive - proof of the fact that entries of the value vector $\mathbf{v}(\beta)$ of $\Gamma_\beta$ lie in $\mathbb{F}$ that is based on a deep result from the field of mathematical logic that is known as {\em Tarski's Principle} \cite{Tarski1951} which states that:

{\em An elmentary sentence which is valid over one real closed field is valid over every real closed field.}

Now, Bewley and Kohlberg \cite{BeKo} showed that Shapley's theorem captured in (\ref{FP}) constitutes a valid elementary sentence over the reals. Since the field $\mathbb{F}$ is known to be real closed, the result follows immediately.

We do not go into all details of the above argument because the main objective of this study, for $\Gamma_\beta$, is to explicitly exhibit the construction of decoupled, bivariate, polynomials whose roots contain the entries of $\mathbf{v}(\beta)$ and to demonstrate the usefulness of Gr{\"o}bner basis techniques in the process. Furthermore, our approach lends itself to extension to the limiting average stochastic games covered in the next section. Since there is no analogue for Shapley's theorem in that case, it is not clear whether Tarski's principle could be used there.

Returning to the main line of discussion, suppose now that the Shapley-Snow kernels $K$ have been fixed and chosen correctly in the sense that $z_s=v_s(\beta); s=1, 2, \ldots, N$ satisfy both (\ref{Redusedfixedpointeq}) and (\ref{Redusedfixedpointeq2ndform}) for $z_0 = \beta$.
That is, $(z_0,\mathbf{z})$, with $\mathbf{z}=[z_1,...,z_N]$, is among the real-valued roots of the system 	of coupled polynomial equations (\ref{syseofcoupledpolyeeqes}) which we now re-write as

	\begin{align}\label{syseofcoupledpolyeeqes2ndver}
		\left\{
 		\begin{array}{c}
 		f_1(z_0,\mathbf{z})=0\\
		\vdots\ \\
		f_N(z_0,\mathbf{z})=0.
 		\end{array}
 		\right.
	\end{align}
	The zero-set of (\ref{syseofcoupledpolyeeqes2ndver}) constitutes an algebraic variety, $W$, which can be difficult to analyse. 	Fortunately, the Gr{\"o}bner basis methods (e.g. see the excellent book of Adams and Loustaunau \cite{lou94}) enable us to, instead,  examine a simpler variety $W_d \supset W$ (of the same dimension as $W$) in a neighbourhood of the solution of the stochastic game that corresponds to the zero-set of suitably constructed decoupled polynomials
	\begin{align}\label{syseofdecoupledpolyeeqes}
	\left\{
 		\begin{array}{c}
 		g_1(z_0,z_1)=0\\
		\vdots\ \\
		g_N(z_0,z_N)=0,
 		\end{array}
 		\right.
	\end{align}
	where each $g_s(z_0, z_s)$ is a bivariate polynomial in only $z_0$ and $z_s$. That is, the zeroes of (\ref{syseofdecoupledpolyeeqes}) contain the zeroes of (\ref{syseofcoupledpolyeeqes2ndver}) and hence the solution of the stochastic game.

	The above simplification follows from the next result which can also be seen as a constructive extension of Lemmata $4.1-4.2$ and Theorem $4.3$ proved in Szczechla et al. \cite{Filar97}. Interestingly, perhaps, the latter results exploit the fact that the contraction property of the Shapley operator $T_\beta$ ensures that the Jacobian of the map corresponding to the fixed point equation is nonsingular.\\

	\begin{theorem}\label{Theorem4.1}
 		Consider a discounted stochastic game $\Gamma_{\beta}$ with all data lying in the ordered field of algebraic numbers. There exist bivariate, decoupled, polynomials $g_s(z_0,z_s)$ in variables $z_0, z_s$, for $s=1,2,\ldots, N$, whose zeroes define the variety $W_d$ containing the solution of the stochastic game, for each $z_0=\beta \in (0,1)$. Furthermore, the corresponding $z_s=v_s(\beta)$ is an algebraic number for each $s=1,2,\ldots, N$.
	\end{theorem}
\medskip
	\begin{proof}
Let $H_\beta:=\{(z_0,z_1,\ldots,z_N) \in \mathbb{C}^{N+1} | z_0=\beta\}$ and consider the variety $W$ that is the zero set of (\ref{syseofcoupledpolyeeqes2ndver}).  From Lemmata $4.1-4.2$ and the proof of Theorem $4.3$ in \cite{Filar97} it follows that $W$ has constant dimension $1$ and that it intersects with $H_\beta$ at a discrete subset of $\mathbb{C}^{N+1}$, that is, for any $\beta \in (0,1)$
\begin{equation}\label{dimzero}
\dim(W\bigcap H_\beta)=0.
\end{equation}

	Now, apply the Buchberger algorithm to find the Gr{\"o}bner basis $GB_1$ of the polynomials $$\lbrace f_1(z_0,\mathbf{z}), \ldots, f_N(z_0,\mathbf{z})\rbrace,$$ with the lexicographic term order $z_1\prec z_2\prec\ldots\prec z_N$.  The last polynomial $g^m(z_0,\mathbf{z}) \in GB_1$ contains the least number of the variables $z_1,\ldots,z_N$ and has coefficients that are polynomial functions of $z_0$.  It is certainly possible to choose some $z_0=\beta^* \in (0,1)$ that is not a root of any of these coefficient polynomials.  But by (\ref{dimzero}) we have that $W_{\beta^*}:=W\bigcap H_{\beta^*}$ is a zero dimensional variety and hence by Corollary $2.2.11$ in \cite{lou94} the polynomial $g^m(\beta^*,\mathbf{z})$ is univariate in $z_1$.  Now, returning to arbitrary $z_0$ in place of $\beta^*,$ define $g_1(z_0, z_1):= g^m(z_0,\mathbf{z})$, a bivariate polynomial in $z_0, z_1$.

Applying  Buchberger algorithm again but with the term order $z_2\prec z_3\prec\ldots\prec z_N\prec z_1$ will, analogously, yield a new Gr{\"o}bner basis $GB_2$ and a bivariate polynomial $g_2(z_0,z_2)$. Continuing in this fashion yields the polynomial 	system (\ref{syseofdecoupledpolyeeqes}).

%	 	The key to finding the polynomials $g_s(z_0,z_s)$, $s=1,2,\ldots, N$ is the finite Buchberger algorithm for 				finding the Gr{\"o}bner basis of the polynomials $$\lbrace f_1(z_0,\mathbf{z}), \ldots, f_N(z_0,\mathbf{z})\rbrace,$$ 				applied $N$ times with different term orders.  Since Gr{\"o}bner basis is non-unique, applying the Buchberger algorithm 				with the term order $z_1\prec z_2\prec\ldots\prec z_N$ yields a Gr{\"o}bner basis the first element of which is 					$g_1(z_0,z_1)$. Applying  that algorithm again but with the term order $z_2\prec z_3\prec\ldots\prec z_N\prec z_1$ will 			yield a new Gr{\"o}bner basis whose first element is $g_2(z_0,z_2)$. Continuing in this fashion yields the polynomial 				system (\ref{syseofdecoupledpolyeeqes}).

		 Note that if $(z_0,\mathbf{z})=(z_0,z_1,\ldots,z_n)$ satisfies (\ref{syseofcoupledpolyeeqes2ndver}), then $g_s(z_0,z_s)=0$, for every $s=1,2,\ldots, N$ since each $g_s$ is part of some Gr{\"o}bner basis of $\lbrace f_1, 					\ldots, f_N \rbrace$. Thus $W\subset W_d$.  In particular, since $(\beta, \mathbf{v}(\beta)) \in W$ by Lemma $4.2$ and Theorem $4.3$ in \cite{Filar97}, it follows that $(\beta, \mathbf{v}(\beta)) \in W_d$ and $g_s(\beta, v_s(\beta))=0,$ for $s=1,\ldots,N$.

		Furthermore, since each $g_s$ is a bivariate polynomial in $z_0, z_s$ with coefficients that are algebraic numbers, it becomes a univariate polynomial in $z_s$ with algebraic coefficients, when evaluated at  $z_0=\beta$ that is algebraic. Thus each $g_s(\beta, z_s)$ is a univariate polynomial with algebraic 					coefficients, hence its roots are also algebraic numbers, since the field of algebraic numbers is algebraically closed. This completes the proof.
	\end{proof}\\

\begin{corollary}\label{Corollary4.1}
Consider a discounted stochastic game $\Gamma_{\beta}$ with all data lying in the ordered field of algebraic numbers $\mathbb{F}$. Then there exists a pair of optimal stationary strategies $(\mu^0, \nu^0)$ with all entries lying in $\mathbb{F}$.
\end{corollary}
\medskip
\begin{proof}
Once we substitute  $\mathbf{u}=\mathbf{v}(\beta)$ with entries lying in $\mathbb{F}$, the Shapley matrix games $R_{\beta}(s,\mathbf{v}(\beta)$ also have all entries that lie in $\mathbb{F}$. The statement of the corollary now follows immediately from the contributions of Weyl \cite{Weyl} and Shapley and Snow \cite{ShapSnow}, see also Proposition $2.1$.
\end{proof}\\

 	\begin{remark}
 	It follows from the finiteness of the Buchberger algorithm (and the fact that there are only $|K|$ kernel selections
 	 $\kappa$), that the correct set of decoupled polynomials (\ref{syseofdecoupledpolyeeqes}) can be found by a finite
 	 algorithm. Of course, the latter could still be of great complexity. Furthermore, the problem of finding the roots of
 	 these polynomials is, in general, not solvable by finite algorithms.
	\end{remark}
%\end{document}
 %**************************************************************************************************************************%
 %**************************************************--- New Section ---******************************************************
 %**************************************************************************************************************************%

\section{Limiting average stochastic games over the ordered field of algebraic numbers}
 	We now consider a stochastic game where the sequence of single stage expected rewards $\lbrace E_{\mu \nu s}(R_t)					\rbrace_{t=0} ^{\infty}$ is aggregated according to
 	\medskip
 	\begin{equation*}
 		v(\mu, \nu, s)= \lim_{\tau\to\infty}\inf\left(\dfrac{1}{\tau +1}\right)\sum_{t=0}^{\tau}E_{\mu \nu s}(R_t).
 	\end{equation*}
 	Such a game is called a {\em limiting average $($or Ces\`{a}ro average$)$ stochastic game} and will be denoted by $\Gamma_{\alpha}$. It is still a zero-sum game to which the minimax solution concept applies. Indeed, the existence of the value vector 	$\mathbf{v}=(v_1,\ldots,v_N)^T$ of $\Gamma_{\alpha}$ was proved by Mertens and Neymann \cite{Mertens}. These authors showed 			that
 	\begin{equation}\label{Mertenseq}
 		\mathbf{v}= \lim_{\beta\to 1^-}\mathbf{v}(\beta),
 	\end{equation}
 	where $\mathbf{v}(\beta)$ is the value vector of the discounted game $\Gamma_{\beta}$. They exploited the following Puiseux series expansion of $\mathbf{v}(\beta)$ that is due to Bewley and Kohlberg \cite{BeKo}
 	\begin{equation}\label{BeKoeq}
 	 \mathbf{v}(\beta)=\sum_{r=0}^{\infty}c_r(s)(1-\beta)^{\frac{r}{M_s}}= c_0(s)+\delta(s,\beta),
 	\end{equation}
 	for each $s=1, 2, \ldots, N$. Here $M_s$ is a natural number and $\delta(s,\beta)\to 0$ as $\beta \to 1$. Clearly, 				(\ref{Mertenseq}) and (\ref{BeKoeq}) imply that the value of the $\Gamma_{\alpha}$ game, starting at state $s$, is given by
 	\begin{equation*}
 	 \mathbf{v}_s= c_0(s),\ s=1, 2, \ldots, N.
 	\end{equation*}
 Unfortunately, we cannot immediately conclude from (\ref{Mertenseq}) that entries of $\mathbf{v}$ lie in the field $\mathbb{F}$ when the data of $\Gamma_\alpha$ lie in $\mathbb{F}$.  This is because a limit of algebraic numbers need not be algebraic. Hence, we needed to establish the next result.\\

 	\begin{theorem}
 	 Consider a limiting average stochastic game $\Gamma_{\alpha}$ with all data lying in the ordered field of algebraic
 	  numbers. Then the value vector $\mathbf{v}$ of $\Gamma_{\alpha}$ has entries $\mathbf{v}_s$ that are also algebraic
 	  numbers, for each $s=1, 2, \ldots, N$.
 	\end{theorem}
 \medskip
 	\begin{proof}
 	 In view of (\ref{Mertenseq}) and (\ref{BeKoeq}), it is sufficient to prove that $c_0(s)$ is an algebraic number for each
 	 $s=1, 2, \ldots, N$. Without loss of generality, we supply a proof for the case of $s=1$ only.

 The proof will be asymptotic,
 	 via a sequence of discounted games $\Gamma_{\beta_{k}}$ such that $\beta_{k}\to 1^{-}$. In particular,
 	 $\lbrace\beta_{k}\rbrace_{k=0}^\infty $ is such that, for each $k $, there is a selection of kernels $\kappa_k\in K $ such,
 	 that the polynomials (\ref{syseofcoupledpolyeeqes2ndver}) define the variety $W^k$ containing $(\beta_k, \mathbf{v}(\beta_{k}))$, where $\mathbf{v}(\beta_{k})$ is the value vector of the game $\Gamma_{\beta_{k}}$. However, since $|K|<\infty $ at least one of these kernels, say $\hat{\kappa},$ must repeat itself infinitely often. Again, without loss of generality, we may assume that the entire sequence  $\lbrace\beta_{k} \rbrace_{k=0}^\infty $ corresponds to the same kernel selection $\hat{\kappa}$ and hence the same variety $\hat{W}$. Now, for the term order $z_1\prec z_2\prec\ldots\prec z_N$ find the Gr{\"o}bner basis of $\lbrace f_1,\ldots, f_N \rbrace,$ and the decoupled polynomial $g_1(z_0,z_2)$, as in the proof of Theorem \ref{Theorem4.1}.  We know that for each $k=0,1,2,\ldots$
 	 \begin{equation}\label{eq5.3}
 	  g_1(\beta_k,v_1(\beta_k))=0.
 	 \end{equation}
 	 Hence $g_1(\beta_k,z_1)$ can be factored as
 	 \begin{equation}\label{eq5.4}
 	  g_1(\beta_k,z_1)=(z_1-v_1(\beta_k))q_1(\beta_k,z_1),
 	 \end{equation}
 	 where $q_1(\beta_k, z_1)$ can be obtained by the long division of $g_1(\beta_k,z_1)$ by the linear factor 		 $(z_1-v_1(\beta_k))$.
 In case $v_1(\beta_k)$ is a root of multiplicity $t >1$, the factor $(z_1-v_1(\beta_k))^{t-1}$ would be included in a factorization of $ g_1(\beta_k,z_1))$ into its irreducible factors. Substituting (\ref{BeKoeq}) into (\ref{eq5.4}) we obtain for each $k=0,1,\ldots$
 	 \begin{equation}\label{eq5.5}
 	  g_1(\beta_k,z_1)=(z_1-c_0(1))q_1(\beta_k,z_1)-\delta(1,\beta_k)q_1(\beta_k,z_1),
 	 \end{equation}
 	 where $q_1(z_0,z_1)$ is a bivariate polynomial in $z_0, z_1$ that has coefficients that are algebraic numbers.

 Next, with $z_0=\beta,$ we can rewrite the latter polynomial as $q_1(\beta,z_1)=(1-\beta)^\ell \bar{q}_1(\beta,z_1)$, where $\bar{q}_1(\beta,z_1)$ does not have a root at $\beta=1$.
 %rational functions of $z_0$, there exists a polynomial $\theta(z_0)$ such that $\bar{q}_1(z_0,z_1)=\theta(z_0)g_1(z_0,z_1)$ is a polynomial in $z_1$ whose coefficients are polynomials in $z_0$. But
 Hence, since $\beta_k \to 1$, we have that $\bar{q}_1(\beta_k,z_1)\to \bar{q}_1(1,z_1)=q^*_1(z_1),$ where the latter is a univariate polynomial in $z_1$ that is free of the parameter $z_0=\beta$.

Now, consider first, the more generic case of $\ell=0$.   Equation (\ref{eq5.5}) becomes
 	 \begin{equation}\label{eq5.6}
 	  g_1(\beta_k,z_1)=(z_1-c_0(1))\bar{q}_1(\beta_k,z_1)-\delta(1,\beta_k)\bar{q}_1(\beta_k,z_1),
 	   \end{equation}
 	   for each $k=0,1,\ldots$ . Passing to the limit with respect to $k$ we obtain
 	   \begin{equation}\label{eq5.7}
 	    \lim_{k\to\infty}\left[g_1(\beta_k,z_1)\right] = g_1(1,z_1) = (z_1-c_0(1))q^*_1(z_1),
 	   \end{equation}
 	   since $\lim_{k\to\infty}\left[\delta(1,\beta_k)\bar{q}_1(\beta_k,z_1)\right] = 0$ and $q^*_1(z_1) \ne 0$.
 	   But by (\ref{eq5.3}) at $z_1=v_1(\beta_k)$, $g_1(\beta_k,v_1(\beta_k))=0$ for each $k$, so that the left side of
 	   (\ref{eq5.7}) is zero. Thus, in the limit
 	   \begin{equation}
 	    0= g_1(1,z_1) = (z_1-c_0(1))q^*_1(z_1)
 	   \end{equation}
 	   and hence $c_0(1)=v_1$  is an algebraic number.

 To complete the proof, we next consider the case where $\ell$ is a positive integer. In this case we have
 \begin{equation}\label{eq5.9}
 	  g_1(\beta_k,z_1)=(z_1-v_1(\beta_k))q_1(\beta_k,z_1)= (z_1-v_1(\beta_k))(1-\beta_k)^\ell \bar{q}_1(\beta_k,z_1),
 	   \end{equation}
 	   for each $k=0,1,\ldots$.  Set $ \bar{g}_1(\beta,z_1):=(z_1-v_1(\beta))\bar{q}_1(\beta,z_1)$, which is a bivariate polynomial in $(\beta, z_1)$ with coefficients that are algebraic numbers. By continuity of polynomials, its limit $\bar{g}_1(1,z_1)$ as $\beta \to 1$ is a polynomial in $z_1$ with coefficients that are algebraic numbers. But, clearly, we still have that  $ 0=\bar{g}_1(\beta_k,v_1(\beta_k))$,
 for each $k=0,1,\ldots$, hence in the limit $ 0=\bar{g}_1(1,v_1)=\bar{g}_1(1,c_0 (1))$ by (\ref{Mertenseq}) and (\ref{BeKoeq}).  Hence, $c_0 (1)$ is a root of a polynomial with algebraic coefficients and is, therefore, an algebraic number.
 	   	\end{proof}

% {\bf Guys, we may include a section (or some remarks only) on the fact that ordered field property over the rational number field $Q$ is very exceptional because in each of the five known classes (3.6) reduce to linear functions and it is hard to imagine the situation where they are higher order polynomials and the solutions would be still rational for all rational $\beta$}.
 	%**************************************************************************************************************************%
 %**************************************************--- New Section ---******************************************************
 %**************************************************************************************************************************%
 \section{Conclusions and examples}

 As hinted in the introduction, we suggest that Theorems $4.1$ and $5.1$ provide a rather comprehensive characterization to the ordered field property problem posed in \cite{Raghavan1991}. This is because the field $\mathbb{F}$ of algebraic numbers is still a countable ordered field that includes rationals $\mathbb{Q}$ and it is unlikely that any smaller extension of rationals would suffice. We have demonstrated that the Gr\"{o}bner basis methods, in principle, allow us to identify exactly polynomials (with integer coefficients) whose roots contain the values of stochastic games $\Gamma_\beta$, for $s=1,\ldots,N$. However, in general, these roots cannot be computed exactly in terms of radicals.  This is illustrated in Example $2$, below, where our decoupled polynomials $g_s( \beta, z_s)$ are quintic. Hence by the classical Abel-Ruffini theorem (e.g., see Theorem $75$ and the following remark in \cite{rotman1998}, page 75), no radical solutions exist for these polynomials.\\

 %For this case, the polynomials turn out to be quintic. Obviously no radical solutions exist due to Abel-Ruffini's theorem (for instance see theorem 75 and the following remark in \cite{rotman1998} page 75). Additional to the non-existance of radical solutions and being able to approximate the solution we need to specify $\beta$. As so, by specifying $\beta=\frac{1}{2}$ the chosen Gr\"{o}bner basis change to:
%
Let us now comment briefly on the structured classes of stochastic games that are known to possess the ordered field property, over the field $\mathbb{Q}$ of rational numbers.  These classes include stochastic games of perfect information originally introduced by Gillette in 1957 \cite{Gillette1957}, separable reward and state independent transition (SER-SIT) games \cite{Sobel1981myopic}, \cite{Parthasarathy1984}, single-controller games \cite{Parthasarathy1981}, switching-controller games \cite{Filar1981} and additive reward and additive transition (ARAT) games \cite{Raghavan85}.  For the sake of completeness, we briefly introduce these classes of games and supply an example of one of them which exposes their simple algebraic structure in the context of Theorem 4.1.\\

 A {\em single controller} stochastic game is a game in which the transition probabilities satisfy the following:
 \begin{align*}
  p(s'|s,a,b_1)=p(s'|s,a,b_2) \; \forall \; s,s'\in \mathbf{S} \; \& \; \forall a\in \mathbf{A}(s), \; b_1,\ b_2\in \mathbf{B}(s).
 \end{align*}
 A {\em switching controller} stochastic game is a generalization of the above for which the state space $S$ can be partitioned into disjoint subsets $S^1$ and $S^2$ such that $\forall \; s' \in S$:
 \begin{align*}
 p(s'|s,a,b_1)=p(s'|s,a,b_2) \; \forall \; s\in S^1 \; \&\ \forall \; a\in \mathbf{A}(s), \; b_1, b_2\in \mathbf{B}(s),\\\nonumber
 p(s'|s,a_1,b)=p(s'|s,a_2,b) \; \forall \; s\in S^2 \; \&\ \forall \; a_1, a_2\in \mathbf{A}(s),\ b\in \mathbf{B}(s).
 \end{align*}
 Incidentally, the {\em perfect information} game is equivalent to a switching controller game in which one of the players has only a single action in each state.

The {\em ARAT} stochastic game  is one in which both the rewards and the transitions can be written as the sum of a term determined by player one and a term determined by player two. That is:
 \begin{align*}
 r(s,a,b)&=r_1(s,a)+r_2(s,b)  \; \forall \; s\in S \; \&\ \forall \; a\in \mathbf{A}(s), \; b\in \mathbf{B}(s),\\\nonumber
 p(s'|s,a,b)&=p_1(s'|s,a)+p_2(s'|s,b)  \; \forall \; s, s' \in S \; \&\ \forall \; a\in \mathbf{A}(s), \; b\in \mathbf{B}(s).
 \end{align*}
Finally, the {\em SER-SIT} stochastic game is one for which the rewards and transition probabilities have the following structure for some scalars $ c(1),c(2),\ldots,c(N)$:
 \begin{align*}
  r(s,a,b)&=c(s)+r(a,b) \; \forall \; s\in S \; \&\ \forall \; a\in \mathbf{A}(s), \; b\in \mathbf{B}(s) \\
  p(s'|s_1,a,b)&=p(s'|s_2,a,b) \; \forall \; s',s_1,s_2\in S \;\&\;\forall\;a\in\mathbf{A}(s)\;b\in\mathbf{B}(s).
  \end{align*}

 Next, we illustrate our approach with two simple examples: the first is a switching controller game and the second is a reward-diagonal game that does not belong to any of the above structured classes.
\newline\\
  {\bf Example 1:}
   We consider the switching controller stochastic game as illustrated in Figure~\ref{fig:Switching}. Its state space $S= S_1 \bigcup S_2 = \{1\} \bigcup \{2\}$.  Each player has three actions in each state, that is,  $\mathbf{A}(1)=\mathbf{A}(2)=\left\lbrace a_1,a_2,a_3\right\rbrace$ and $\mathbf{B}(1)=\mathbf{B}(2)=\left\lbrace b_1,b_2,b_3\right\rbrace$. For example, in Figure(\ref{fig:Switchingsub1}) if players one and two choose the first column and row, respectively, the immediate reward will be $-2$ for player one and consequently $+2$ for player two. Furthermore, at the next stage the game stays in state $1$ with probability $\frac{3}{10}$ and moves to state $2$ with probability $\frac{7}{10}$.

    We solved this game in MAPLE symbolic computing environment.  In particular, postponing the normalization by the factor $(1-\beta)$, Shapley's auxiliary matrix games were evaluated and equations (\ref{Redusedfixedpointeq}) were manipulated to arrive at the system of two polynomial equations of the form (\ref{Redusedfixedpointeq2ndform}). The latter are (with $z_0=\beta$):
\begin{align*}
f_1(\beta, z_1, z_2) =&  32-12 z_1+\frac{22}{5} \beta z_1+\frac{38}{5} \beta z_2\\
f_2(\beta, z_1, z_2) =& 9 z_2-\frac{81}{20} \beta z_1-\frac{99}{20} \beta z_2.
\end{align*}
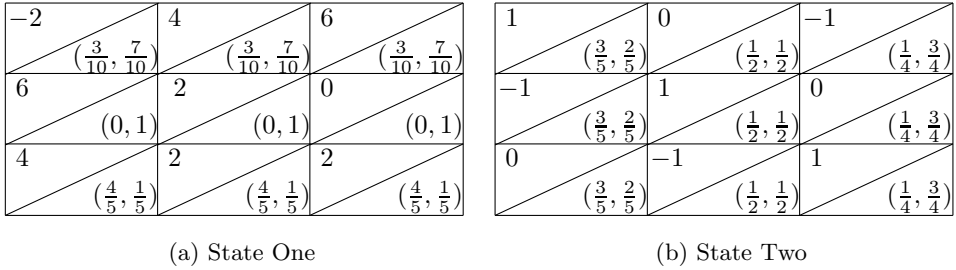
\begin{figure}[h!]
\centering
\begin{subfigure}{0.5\textwidth}
  \centering
\begin{tikzpicture}%[thick]
  \matrix (mat) [%
    matrix of nodes,
    nodes in empty cells,
    text width=1.79cm,
    text height=5pt,
    text depth=15pt,
    text badly centered
  ](mat)
  {%
     \llap{\raisebox{0pt}{$-2$}\hspace{15pt}}\llap{\raisebox{-14.4pt}{$(\frac{3}{10},\frac{7}{10})$}\hspace{-29.4pt}}&
      \llap{\raisebox{0pt}{$4$}\hspace{20pt}}\llap{\raisebox{-14.4pt}{$(\frac{3}{10},\frac{7}{10})$}\hspace{-29.4pt}}& \llap{\raisebox{0pt}{$6$}\hspace{20pt}}\llap{\raisebox{-14.4pt}{$(\frac{3}{10},\frac{7}{10})$}\hspace{-29.4pt}}& \\
    \llap{\raisebox{0pt}{$6$}\hspace{20pt}}\llap{\raisebox{-14pt}{$(0,1)$}\hspace{-29pt}}&
      \llap{\raisebox{0pt}{$2$}\hspace{18pt}}\llap{\raisebox{-14pt}{$(0,1)$}\hspace{-29pt}}& \llap{\raisebox{0pt}{$0$}\hspace{20pt}}\llap{\raisebox{-14pt}{$(0,1)$}\hspace{-29pt}}& \\
    \llap{\raisebox{0pt}{$4$}\hspace{20pt}}\llap{\raisebox{-14pt}{$(\frac{4}{5},\frac{1}{5})$}\hspace{-29pt}}&
      \llap{\raisebox{0pt}{$2$}\hspace{20pt}}\llap{\raisebox{-14pt}{$(\frac{4}{5},\frac{1}{5})$}\hspace{-29pt}}& \llap{\raisebox{0pt}{$2$}\hspace{20pt}}\llap{\raisebox{-14pt}{$(\frac{4}{5},\frac{1}{5})$}\hspace{-29pt}}& \\
  };
  % horizontal lines
  \foreach \i in {1,2,3}
        \draw (mat-\i-1.north west) -- (mat-\i-3.north east);
  \draw (mat-3-1.south west) -- (mat-3-3.south east);
  % vertical lines
  \foreach \j in {1,2,3}
    \draw (mat-1-\j.north west) -- (mat-3-\j.south west);
  \draw (mat-1-3.north east) -- (mat-3-3.south east);
  % diagonal line
  \draw (mat-1-1.north east) -- (mat-1-1.south west);
  \draw (mat-1-2.north east) -- (mat-1-2.south west);
  \draw (mat-1-3.north east) -- (mat-1-3.south west);
  \draw (mat-2-1.north east) -- (mat-2-1.south west);
  \draw (mat-2-2.north east) -- (mat-2-2.south west);
  \draw (mat-2-3.north east) -- (mat-2-3.south west);
  \draw (mat-3-1.north east) -- (mat-3-1.south west);
  \draw (mat-3-2.north east) -- (mat-3-2.south west);
  \draw (mat-3-3.north east) -- (mat-3-3.south west);
\end{tikzpicture}
\centering
  \caption{State One}
  \label{fig:Switchingsub1}
\end{subfigure}%
\begin{subfigure}{.5\textwidth}
  \centering
 \begin{tikzpicture}%[thick]
  \matrix (mat) [%
     matrix of nodes,
    nodes in empty cells,
    text width=1.79cm,
    text height=5pt,
    text depth=15pt,
    text badly centered
  ]
  {%
     \llap{\raisebox{0pt}{$1$}\hspace{20pt}}\llap{\raisebox{-14pt}{$(\frac{3}{5},\frac{2}{5})$}\hspace{-29pt}}&
      \llap{\raisebox{0pt}{$0$}\hspace{20pt}}\llap{\raisebox{-14pt}{$(\frac{1}{2},\frac{1}{2})$}\hspace{-29pt}}& \llap{\raisebox{0pt}{$-1$}\hspace{15pt}}\llap{\raisebox{-14pt}{$(\frac{1}{4},\frac{3}{4})$}\hspace{-29pt}}& \\
    \llap{\raisebox{0pt}{$-1$}\hspace{15pt}}\llap{\raisebox{-14pt}{$(\frac{3}{5},\frac{2}{5})$}\hspace{-29pt}}&
      \llap{\raisebox{0pt}{$1$}\hspace{20pt}}\llap{\raisebox{-14pt}{$(\frac{1}{2},\frac{1}{2})$}\hspace{-29pt}}& \llap{\raisebox{0pt}{$0$}\hspace{20pt}}\llap{\raisebox{-14pt}{$(\frac{1}{4},\frac{3}{4})$}\hspace{-29pt}}& \\
    \llap{\raisebox{0pt}{$0$}\hspace{20pt}}\llap{\raisebox{-14pt}{$(\frac{3}{5},\frac{2}{5})$}\hspace{-29pt}}&
      \llap{\raisebox{0pt}{$-1$}\hspace{15pt}}\llap{\raisebox{-14pt}{$(\frac{1}{2},\frac{1}{2})$}\hspace{-29pt}}& \llap{\raisebox{0pt}{$1$}\hspace{20pt}}\llap{\raisebox{-14pt}{$(\frac{1}{4},\frac{3}{4})$}\hspace{-29pt}}& \\
  };
  % horizontal lines
  \foreach \i in {1,2,3}
        \draw (mat-\i-1.north west) -- (mat-\i-3.north east);
  \draw (mat-3-1.south west) -- (mat-3-3.south east);
  % vertical lines
  \foreach \j in {1,2,3}
    \draw (mat-1-\j.north west) -- (mat-3-\j.south west);
  \draw (mat-1-3.north east) -- (mat-3-3.south east);
  % diagonal line
  \draw (mat-1-1.north east) -- (mat-1-1.south west);
  \draw (mat-1-2.north east) -- (mat-1-2.south west);
  \draw (mat-1-3.north east) -- (mat-1-3.south west);
  \draw (mat-2-1.north east) -- (mat-2-1.south west);
  \draw (mat-2-2.north east) -- (mat-2-2.south west);
  \draw (mat-2-3.north east) -- (mat-2-3.south west);
  \draw (mat-3-1.north east) -- (mat-3-1.south west);
  \draw (mat-3-2.north east) -- (mat-3-2.south west);
  \draw (mat-3-3.north east) -- (mat-3-3.south west);
\end{tikzpicture}
\centering
  \caption{State Two}
  \label{fig:Switchingsub2}
\end{subfigure}
\caption{A switching controller stochastic game.}
\label{fig:Switching}
\end{figure}

To compute the  Gr\"{o}bner  bases, we first used the lexicographic term order $z_1\prec z_2$ and then applied the same process again but with the term order $z_2\prec z_1$. From the two  Gr\"{o}bner  bases $GB_1, \; GB_2$ so obtained, we selected the two bivariate polynomials:
\begin{itemize}
\item For $z_1\prec z_2$:
\begin{equation}
 g_1(\beta, z_1)=5\beta^2 z_1+55\beta z_1-88\beta-60 z_1+160
\end{equation}
%\begin{equation}
% g_2(\beta, z_2)=5\beta^2 z_2+55 \beta z_2+72 \beta-60 z_2
%\end{equation}
\item For $z_2\prec z_1$:
\begin{equation}
 g_2(\beta, z_2)=5\beta^2 z_2+55 \beta z_2+72 \beta-60 z_2.
\end{equation}
%\begin{equation}
% g_1(\beta, z_1, z_2)=5\beta^2 z_1+55\beta z_1-88\beta-60 z_1+160
%\end{equation}
\end{itemize}
It is not a coincidence, that the latter are linear functions in $z_1$ and $z_2$, respectively. Evidently, this is a common feature for all the five structured classes of stochastic games listed above.
The zeroes of the above polynomials occur when
\begin{align}
z_1=&\frac{8}{5}\frac{11\beta-20}{(\beta+12)(\beta-1)},\\
z_2=&-\frac{72}{5}\frac{\beta}{(\beta+12)(\beta-1)}.
\end{align}
which are clearly positive for $\beta\in\left[ 0,1\right)$. They are also rational for rational $\beta$ and algebraic for algebraic $\beta$.
Now, for $\beta \in (0,1)$ the entries of the value vector $\mathbf{v}_\beta$ of $\Gamma_\beta$ are obtained by normalization by the factor $(1-\beta)$, namely
\begin{align}
v_1(\beta)=(1-\beta)z_1=&\frac{160-88\beta}{5(\beta+12)},\\
v_2(\beta)=(1-\beta)z_2=&\frac{72\beta}{5(\beta+12)}.
\end{align}

{\bf Example 2:}
The next example is related to the same approach for a reward-diagonal stochastic game with data as given
in Figure~\ref{fig:MIXED}.
\begin{figure}[h!]
\centering
\begin{subfigure}{0.5\textwidth}
  \centering
\begin{tikzpicture}%[thick]
  \matrix (mat) [%
    matrix of nodes,
    nodes in empty cells,
    text width=1.79cm,
    text height=5pt,
    text depth=15pt,
    text badly centered
  ]
  {%
     \llap{\raisebox{0pt}{$18$}\hspace{18pt}}\llap{\raisebox{-15pt}{$(0,1)$}\hspace{-29pt}}&
      \llap{\raisebox{0pt}{$0$}\hspace{20pt}}\llap{\raisebox{-14pt}{$(\frac{1}{5},\frac{4}{5})$}\hspace{-29pt}}& \llap{\raisebox{0pt}{$0$}\hspace{20pt}}\llap{\raisebox{-15pt}{$(1,0)$}\hspace{-29pt}}& \\
    \llap{\raisebox{0pt}{$0$}\hspace{20pt}}\llap{\raisebox{-14pt}{$(\frac{2}{5},\frac{3}{5})$}\hspace{-29pt}}&
      \llap{\raisebox{0pt}{$12$}\hspace{18pt}}\llap{\raisebox{-14pt}{$(\frac{1}{3},\frac{2}{3})$}\hspace{-29pt}}& \llap{\raisebox{0pt}{$0$}\hspace{20pt}}\llap{\raisebox{-14pt}{$(\frac{3}{5},\frac{2}{5})$}\hspace{-29pt}}& \\
    \llap{\raisebox{0pt}{$0$}\hspace{20pt}}\llap{\raisebox{-14pt}{$(\frac{3}{5},\frac{2}{5})$}\hspace{-29pt}}&
      \llap{\raisebox{0pt}{$0$}\hspace{20pt}}\llap{\raisebox{-14pt}{$(\frac{2}{5},\frac{3}{5})$}\hspace{-29pt}}& \llap{\raisebox{0pt}{$6$}\hspace{20pt}}\llap{\raisebox{-14pt}{$(0,1)$}\hspace{-29pt}}& \\
  };
  % horizontal lines
  \foreach \i in {1,2,3}
        \draw (mat-\i-1.north west) -- (mat-\i-3.north east);
  \draw (mat-3-1.south west) -- (mat-3-3.south east);
  % vertical lines
  \foreach \j in {1,2,3}
    \draw (mat-1-\j.north west) -- (mat-3-\j.south west);
  \draw (mat-1-3.north east) -- (mat-3-3.south east);
  % diagonal line
  \draw (mat-1-1.north east) -- (mat-1-1.south west);
  \draw (mat-1-2.north east) -- (mat-1-2.south west);
  \draw (mat-1-3.north east) -- (mat-1-3.south west);
  \draw (mat-2-1.north east) -- (mat-2-1.south west);
  \draw (mat-2-2.north east) -- (mat-2-2.south west);
  \draw (mat-2-3.north east) -- (mat-2-3.south west);
  \draw (mat-3-1.north east) -- (mat-3-1.south west);
  \draw (mat-3-2.north east) -- (mat-3-2.south west);
  \draw (mat-3-3.north east) -- (mat-3-3.south west);
\end{tikzpicture}
\centering
  \caption{State One}
  \label{fig:MIXEDsub1}
\end{subfigure}%
\begin{subfigure}{.5\textwidth}
  \centering
 \begin{tikzpicture}%[thick]
  \matrix (mat) [%
     matrix of nodes,
    nodes in empty cells,
    text width=1.79cm,
    text height=5pt,
    text depth=15pt,
    text badly centered
  ]
  {%
     \llap{\raisebox{0pt}{$3$}\hspace{20pt}}\llap{\raisebox{-14pt}{$(\frac{1}{5},\frac{4}{5})$}\hspace{-29pt}}&
      \llap{\raisebox{0pt}{$0$}\hspace{20pt}}\llap{\raisebox{-14.4pt}{$(\frac{7}{10},\frac{3}{10})$}\hspace{-29.5pt}}& \llap{\raisebox{0pt}{$0$}\hspace{20pt}}\llap{\raisebox{-14.4pt}{$(\frac{3}{10},\frac{7}{10})$}\hspace{-29.5pt}}& \\
    \llap{\raisebox{0pt}{$0$}\hspace{20pt}}\llap{\raisebox{-14pt}{$(\frac{1}{2},\frac{1}{2})$}\hspace{-29pt}}&
      \llap{\raisebox{0pt}{$2$}\hspace{20pt}}\llap{\raisebox{-14pt}{$(\frac{2}{5},\frac{3}{5})$}\hspace{-29pt}}& \llap{\raisebox{0pt}{$0$}\hspace{20pt}}\llap{\raisebox{-14pt}{$(\frac{2}{5},\frac{3}{5})$}\hspace{-29pt}}& \\
    \llap{\raisebox{0pt}{$0$}\hspace{20pt}}\llap{\raisebox{-14.4pt}{$(\frac{3}{10},\frac{7}{10})$}\hspace{-29.5pt}}&
      \llap{\raisebox{0pt}{$0$}\hspace{20pt}}\llap{\raisebox{-14pt}{$(\frac{1}{5},\frac{4}{5})$}\hspace{-29pt}}& \llap{\raisebox{0pt}{$1$}\hspace{20pt}}\llap{\raisebox{-15pt}{$(0,1)$}\hspace{-29pt}}& \\
  };
  % horizontal lines
  \foreach \i in {1,2,3}
        \draw (mat-\i-1.north west) -- (mat-\i-3.north east);
  \draw (mat-3-1.south west) -- (mat-3-3.south east);
  % vertical lines
  \foreach \j in {1,2,3}
    \draw (mat-1-\j.north west) -- (mat-3-\j.south west);
  \draw (mat-1-3.north east) -- (mat-3-3.south east);
  % diagonal line
  \draw (mat-1-1.north east) -- (mat-1-1.south west);
  \draw (mat-1-2.north east) -- (mat-1-2.south west);
  \draw (mat-1-3.north east) -- (mat-1-3.south west);
  \draw (mat-2-1.north east) -- (mat-2-1.south west);
  \draw (mat-2-2.north east) -- (mat-2-2.south west);
  \draw (mat-2-3.north east) -- (mat-2-3.south west);
  \draw (mat-3-1.north east) -- (mat-3-1.south west);
  \draw (mat-3-2.north east) -- (mat-3-2.south west);
  \draw (mat-3-3.north east) -- (mat-3-3.south west);
\end{tikzpicture}
\centering
  \caption{State Two}
  \label{fig:MIXEDsub2}
\end{subfigure}
\caption{A reward-diagonal stochastic game.}
\label{fig:MIXED}
\end{figure}
By following the same method as in previous example we derived Shapley's polynomials (\ref{Redusedfixedpointeq2ndform}). The latter are (with $z_0=\beta$):
\begin{align*}
f_1(\beta, z_1, z_2)=&-\frac{2}{75} \beta^2 z_1^3+\frac{4}{75} \beta^2 z_1^2 z_2-\frac{89}{5} \beta z_1^2-\frac{2}{75} z_1 \beta^2 z_2^2+\frac{89}{5} z_1 \beta z_2+396 z_1+\frac{1}{375} \beta^3 z_1^2 z_2+9 \beta^2 z_1^2\\
&-\frac{11}{375} \beta^3 z_1 z_2^2-\frac{1}{5} \beta^2 z_1 z_2-144 \beta z_1+\frac{7}{375} \beta^3 z_2^3-\frac{44}{5} \beta^2 z_2^2-252 \beta z_2-1296+\frac{1}{125} \beta^3 z_1^3
 \\
f_2(\beta, z_1, z_2)=& -\frac{6}{25} z_1 \beta^2 z_2^2-2 z_1 \beta z_2+\frac{3}{25} \beta^2 z_2^3+2 \beta z_2^2+11 z_2+\frac{3}{25} \beta^2 z_1^2 z_2+\frac{33}{500} \beta^3 z_1^2 z_2+\frac{69}{100} \beta^2 z_1^2\\
&+\frac{27}{500} \beta^3 z_1 z_2^2+\frac{31}{50} \beta^2 z_1 z_2-\frac{8}{5} \beta z_1-\frac{31}{500} \beta^3 z_1^3-\frac{29}{500} \beta^3 z_2^3-\frac{131}{100} \beta^2 z_2^2-\frac{47}{5} \beta z_2-6.
\end{align*}
%To compute the  Gr\"{o}bner basis, again we first use the term order $z_1\prec z_2$ and then we apply the same process with the term order $z_2\prec z_1$. From the two lists that we obtain in this way, we pick bivariate polynomials in $\beta, z_1$ and $\beta, z_2$ respectively:
To compute the  Gr\"{o}bner  bases, we first used the lexicographic term order $z_1\prec z_2$ and then applied the same process again but with the term order $z_2\prec z_1$. From the two  Gr\"{o}bner  bases $GB_1, \; GB_2$ so obtained, we selected the two bivariate polynomials:
\begin{itemize}
\item For $z_1\prec z_2$:
\begin{align*}
g_1(\beta,z_1)=&\left( 573788072 \beta^8-220687720 \beta^7-7503382480 \beta^6+20744645680 \beta^5\right. \\
&\left.-23613586040 \beta^4+12667475128 \beta^3-2648252640 \beta^2\right)z_1^5+\left( 17339824845 \beta^7\right. \\
&\left.+777484096540 \beta^6-5051043271810 \beta^5+12082535624940 \beta^4\right. \\
&\left.-13976285582035 \beta^3+7917677944720 \beta^2-1767708637200 \beta\right)  z_1^4\\
&+\left( -132203187917 \beta^6-4061892203389 \beta^5+41976312974469 \beta^4\right. \\
&\left.-139562105527903 \beta^3+209430408589540 \beta^2-146977072348800 \beta\right. \\
&\left.+39326551704000\right)  z_1^3+\left( 4540064732898 \beta^5-32189023528776 \beta^4\right. \\
&\left.+111423218218458 \beta^3-210428562654180 \beta^2+192643241103600 \beta\right. \\
&\left.-65988937872000\right)  z_1^2+\left( -36494410557024 \beta^4+197044804007424 \beta^3\right. \\
&\left.-329307266563200 \beta^2+153965516032800 \beta+14791357080000\right)  z_1\\
&+89067118187808 \beta^3-570684050648640 \beta^2+1146714437289600 \beta \\
&-720144590112000.
\end{align*}
\item For $z_2\prec z_1$:
\begin{align*}
g_2(\beta,z_2)=&\left( 573788072 \beta^8-220687720 \beta^7-7503382480 \beta^6+20744645680 \beta^5-23613586040 \beta^4\right. \\
&\left. +12667475128 \beta^3-2648252640 \beta^2\right) z_2^5+\left( -204848895 \beta^7-858363627960 \beta^6\right. \\
&+3391365456790 \beta^5-4977728569660 \beta^4+3171701981265 \beta^3-682632847540 \beta^2\\
&\left.-44137544000 \beta\right)z_2^4+\left( -428935189842 \beta^6-27966474140289 \beta^5+100992652545319 \beta^4\right. \\
&\left.-130647004403803 \beta^3+72083867190015 \beta^2-13791349509400 \beta-242756492000\right)z_2^3 \\
&+\left( -11840380819272 \beta^5-294095993630166 \beta^4+1011540632651148 \beta^3\right. \\
&\left.-1175534529822210 \beta^2+552033040308000 \beta-82102768687500\right)  z_2^2\\
&+\left( -150578129579184 \beta^4-965930452772256 \beta^3+3604616402786640 \beta^2 \right. \\
&\left. -3768621176239200 \beta+1280513355804000\right)  z_2 -737655098899392 \beta^3\\
&+906896369212560 \beta^2+449706809325600 \beta-673995164922000.
\end{align*}
\end{itemize}
For this case, the polynomials turn out to be quintic. As mentioned earlier, no radical solutions exist due to Abel-Ruffini's theorem. In order to find approximate real roots of these polynomials, we need to specify $\beta$. For instance, by setting $\beta=\frac{1}{2}$ the above polynomials simplify to:
\begin{align*}
g_1(\frac{1}{2}, z_1)=&-\frac{5359807932115}{128}  z_1^4-\frac{733786669}{32}  z_1^5+\frac{207671951406997}{64}  z_1^3\\
&+31796998296714  z_1-278324994355884-\frac{163463882465331}{16}  z_1^2\\
g_2(\frac{1}{2},z_2)=&-\frac{733786669}{32}  z_2^5+\frac{61773492718827}{8}  z_2^2-\frac{1895121071975}{128}  z_2^4\\
&+167204428685829  z_2-314624555318484-\frac{33907289629}{2}  z_2^3.
%g_2(z_1,\frac{1}{2})=&-\frac{5359807932115}{128}  z_1^4-\frac{733786669}{32}  z_1^5+\frac{207671951406997}{64}  z_1^3\\
%&+31796998296714  z_1-278324994355884-\frac{163463882465331}{16}  z_1^2
\end{align*}
Approximating their real roots in MAPLE we obtain:\\
$$z_1=-1900.653702, 4.969443147, 71.41363761$$
and \\
$$z_2=-643.7311436, 1.742768501, 28.64701004.$$
Note that the unique solution of the stochastic game $\Gamma_{\frac{1}{2}}$ is obtained only with the help of one particular pair of the above real roots. Indeed, it follows that
\begin{align}
v_1(\frac{1}{2})=(1-\frac{1}{2})z_1=&\frac{1}{2}4.969443147=2.484721573\\
v_2(\frac{1}{2})=(1-\frac{1}{2})z_2=&\frac{1}{2}1.742768501=0.67138250.
\end{align}
Formal verification that the above are indeed values of $\Gamma_{\frac{1}{2}}$, involves substituting these into Shapley's auxilliary matrix games and finding their values.

\newpage
\bibliographystyle{siam}
\bibliography{OrderedFieldDec14}
\end{document}